\documentclass[12pt]{article}

\usepackage[utf8]{inputenc}
\usepackage[margin=3cm]{geometry}
\usepackage{amsfonts, amsmath, amssymb, amsthm}
\usepackage{indentfirst}

\newtheorem{Thm}{Theorem}

\newtheorem{Lem}{Lemma}

\newenvironment{Prf}{\par\noindent{\it Proof:}\rm}{\newline\rightline{\textbf{QED}}\par}
\newtheorem{Cor}{Corollary}

\newcommand{\fin}{[\omega]^{<\omega}}

\newcommand{\unc}{\rm{unc}}
\newcommand{\Span}{\rm{Span}}
\newcommand{\Cl}{\rm{Cl}}

\title{\large A simple, combinatorial proof of Gowers Dichotomy for Banach Spaces \\ Une preuve combinatoire simple de la dichotomie de Gowers pour les espaces de Banach}

\author{\large Ryszard Frankiewicz, Sławomir Kusiński}

\date{}

\begin{document}
  
  \maketitle
  
  \begin{abstract}
  	In this paper we present a simple proof of Gowers Dichotomy which states that every infinite dimensional Banach Space has a subspace which either contains an unconditional basic sequence or is hereditarily indecomposable. Our approach is purely combinatorial and mainly based on work of Ellentuck, Galvin and Prikry in infinite Ramsey theory.
      
    \begin{center}
        \textbf{Résumé}
    \end{center}
    Dans cet article, nous présentons une preuve simple du théorème de dichotomie de Gowers qui stipule que tout espace de Banach de dimension infinie possède un sous-espace qui est ou bien à base inconditionnelle ou bien héréditairement indécomposable.
    Notre approche est purement combinatoire et principalement fondée sur les travaux d'Ellentuck, Galvin et Prikryin en théorie de Ramsey infini.
  \end{abstract}

  \section{Introduction}
  
    Let \( X \) be an infinite dimensional Banach space. We will write \( Y < X \) when \( Y \) is an infinite dimensional, closed subspace of \( X \). Our goal is to give a simple proof and show a combinatorial nature of the following theorem due to Gowers\cite{Gow}:
    \begin{Thm}
      	There exists \( Y < X \) such that either \( Y \) has an unconditional basis or \( Y \) is hereditarily indecomposable.
    \end{Thm}
    Deep connections of this theorem with combinatorics and in particular Ramsey theory have been known from the beginning. Initially our objective was to present our argument in close resemblance to the work of Ellentuck\cite{Ell}, but along the way we needed to introduce some functional analysis into our combinatorics. The whole paper is expressed in combinatorial language in order to emphasise aforementioned connections. It is worth mentioning that the method presented does not generalize directly to non-separable spaces. Such a generalization is a subject for further investigation. 
    
    By \( \omega \) we will denote the set of non-negative integers. For any infinite set \(A\) by \( [A]^{\omega} \) we will denote the family of all infinite subsets of \(A\) and by \( [A]^{<\omega} \) we will denote the family of all finite subsets of \( A \). Recall that a Banach space is hereditarily indecomposable if none of its infinite dimensional subspaces is isomorphic to a topological direct sum of two infinite dimensional subspaces. A pair \( (U, V) \) of subspaces of \( X \) is \( C \)-uncondtional if
    \[
      \forall_{u \in U, v \in V} \| u - v \| \le C\| u + v \|.
    \]
    We can now define \(\unc(U,V) = \inf (\{ C \colon \mbox{ \((U,V)\) is \(C\)-unconditional } \} \cup \{ \infty \}) \). \( X \) is hereditarily indecomposable iff for all \( U, V < X \) we have \( \unc(U,V) = \infty \). We can think of uncoditionality \( \unc \) as of the measure of closeness between space. It is greater the closer the spaces are together and it is always infinite if they intersect on a non-trivial space. 
    
    A sequence \((x_n)_{n<\omega}\) is called unconditional if there exists \(C\) such that for all \( \sum\limits_{n = 1}^\infty a_n x_n \) and all \( \varepsilon_n \in \{-1, 1\} \) we have
    \[
      \| \sum\limits_{n = 1}^\infty \varepsilon_n a_n x_n \| \le C\| \sum\limits_{n = 1}^\infty a_n x_n \|.
    \]
    An unconditional basis is a Schauder basis that is unconditional.
    
    Because of the existential quantifier in the theorem we can without loss of generality assume \( X \) to be separable. Then there exists \( Q \subseteq X \) such that \( \Cl(Q) = S_X(0,1) \) and \( Q \) is countable. Let \( Q = \{ x_i \colon i \in \omega \} \). 
    
    Then for \( A \subset \omega \) we can define
    \[
      Q_A = \{ x_n \colon n \in A \}
    \] and
    \[
      X_A = \Cl(\Span(A_A)).
    \]
    With no loss of generality we may assume that for any \( A \subseteq \omega \) such that \( |A| > 1 \) the closure \( \Cl(Q \cap X_A) = S_{X_A}(0,1) \).
    
    Our proof methodology is combinatorial, heavily inspired by works of Ellentuck\cite{Ell} and Galvin-Prikry\cite{GP}. For \( s, A \subseteq \omega \) such that \( s \subseteq A \), \( s \in \fin \) and \( X_A \) is infinitely dimensional we define
    \[
        (s, A)^\omega = \{ S \subseteq A \colon s \subseteq S, \: \max(s) < \min(S \setminus s), \: X_S < X \}
    \]
    that is the set all infinite subsets of \( A \) generating an infinitely dimensional subspace and having \( s \) as their initial segment. It is also worth noting that the sets \( (s, A)^\omega \) form a basis for the topology on \( 2^\omega \) stronger than the product topology.
    For \( S, T \in (\emptyset, A)^\omega \) we can define
    \[
      \unc(S, T) = \unc(X_S, X_T).
    \]
    For any \( A \subseteq \omega \) (finite or infinite alike) we can define
    \[
        \tilde{A} = \{ i \in \omega: x_i \in X_A \}.
    \]
    Observe that \( |\tilde{A}| = \omega \) iff \( |A| > 1 \). We can view this operation as the completion of a set \(A\) to a maximal dense subset of the sphere \( S_{X_A}(0,1) \). As we will see later it will be vital to our considerations.

    For earlier approaches to the topic the reader might want to see \cite{RF}, where countable linear spaces over rational numbers were used, as well as \cite{Maur}.
	  
  \section{The positive case}
  
    Let \( C \in (1;+\infty), \: A \in (\emptyset, \omega)^\dag, \: s, t \in [A]^{<\omega} \). We will say that \( A \) is \(C\)-accepting \( (s, t) \) if
    \begin{equation}
      \forall_{B \in (s \cup t, \tilde{A})^\omega} \exists_{S \in (s, \tilde{B})^\omega, T \in (t, \tilde{B})^\omega}
      \unc(S, T) < C.
    \end{equation}
    
    That means the finite dimensional spaces are accepted if they can be extended to infinite dimensional subspaces which closeness is less than \(C\).
    
    We will say that \( A \) is \(C\)-rejecting \( (s, t) \) if
    \begin{equation}
        \forall_{S \in (s, \tilde{A})^\omega, T \in (t, \tilde{A})^\omega}
        \unc(S, T) \ge C.
    \end{equation}
    Observe that for fixed \( A \) and \( C \) there very well may exist \( (s,t) \) that is neither accepted nor rejected, ie there is some grey space between those two notions. \( A \) will be called \( C \)-dichotomous if it either accepts or rejects each pair of its finite subsets.
    
    Now let \( s \in [A]^{<\omega} \). We will say that \( A \) is \(C\)-accepting \( s \) if
    \begin{equation}
      \forall_{t \in P(s)}  \mbox{ \(A\) is \(C\)-accepting } (t, s \setminus t).
    \end{equation}
    
     One of our main objectives will be showing that the emptyset being accepted by some infinite dimensional subspace is sufficient for the existence of an unconditional basis.
     
     \begin{Lem}
         If \( A \in (\emptyset, \omega)^\omega \), \( A = \tilde{A} \) and \( C > 1 \) then there exists a set \( A^* \in (\emptyset, A)^\omega \) such that for any pair \( (s,t) \) of finite subsets of \(A^*\), it is either \( C \)-accepting or \( C \)-rejecting \( (s,t) \), ie \( A^* \) is \( C \)-dichotomous.
     \end{Lem}
     \begin{Prf}
         As the set \( [A]^{<\omega} \times [A]^{<\omega} \) is countable we can enumerate it \( (s_n, t_n)_{n \in \omega} \). Let \( A_0 = \tilde{A} \) and assume that the sets \( A_0 \supseteq \ldots \supseteq A_n \) are defined.
         
         If \( A_n \) is \(C\)-accepting \( (s_n, t_n) \) then take \( A_{n+1} = A_n \). In other case we have
         \[
            \exists_{B \in (\emptyset, A_n)^\omega} \forall_{S \in (s, \tilde{B})^\omega} \forall_{T \in (t, \tilde{B})^\omega} \unc(S, T) \ge C.
         \]
         It follows that \( B \) is \(C\)-rejecting \( (s_n, t_n) \) and so is \( \tilde{B} \). As \( A_n = \tilde{A_n} \) we can take \( A_{n+1} = \tilde{B} \).
         
         With the sets \( A_n \) defined for \( n \in \omega \) if the sequence \( (A_n)_{n \in \omega} \) is constant for \( n \) greater than some \( N \) then we can take \( A^* = A_{N+1} \) and it will accept all but finitely many pairs of its finite subsets. Otherwise let us pick \( i_n \in A_n \) in such a way that the set \( B = \{ i_n \colon n \in \omega \} \) is linearly independent. Clearly \( B \) is almost contained in each \( A_n \), ie \( |B \ A_n| < \omega \). Then we can take \( A^* = \tilde{B} \). Of course \( A^* \in (\emptyset, A)^\omega \) as \( B \subseteq A \) and \( A = \tilde{A} \).
         Now if \( s,t \in [ A^* ]^{<\omega} \subseteq  [ A^* ]^{<\omega} \) then \( (s,t) = (s_n, t_n) \) for some \( n \in \omega \). \(A^*\) is \(C\)-accepting \( (s_n, t_n) \) if \( A_n \) is \(C\)-accepting \( (s_n, t_n) \) and otherwise \( A^* \) is \(C\)-rejecting \( (s_n, t_n) \) as \( A_{n+1} \) is rejecting it. 
     \end{Prf}
 
    From now on let us assume that \( A = \tilde{A} \) and \( A \) is always either \( C \)-accepting or \( C \)-rejecting each pair of its finite subsets.
    
    \begin{Lem}
      If \(A\) is \(C\)-accepting \((s, t)\) then
      \begin{equation*}
        \forall_{B \in (s \cup t, \tilde{A})^\omega}
        \exists_{B^* \in (\emptyset, \tilde{B})^\omega} \bigg(\max(s \cup t) < \min(B^*) \wedge \forall_{i \in B^*}
        \mbox{ \(A\) is \(C\)-accepting \((s \cup \{i\}, t)\) }
        \bigg)
      \end{equation*}
    \end{Lem}
    \begin{Prf}
      Fix \( B \in (s \cup t, A)^\omega \). Without loss of generality we might assume \( B = \tilde{B} \). As \(A\) is \(C\)-accepting \((s,t)\) then for \(B\) there do exist \(S \in (s, \tilde{A})^\omega\) and \(T \in (t, \tilde{A})^\omega \) such that \( \unc(S,T) \). That means that for any \( i \in S \) such that \( i > \max(s \cup t) \) the set \( A \) cannot be \(C\)-rejecting \( (s \cup \{i\}, t) \) and it thus must be \( C \)-accepting it.
        
      Then we can take \( B^* = S \setminus \max(s \cup t) \).
    \end{Prf}

    \begin{Lem}
      If \(A\) is \(C\)-accepting \(s\) then
      \begin{equation*}
        \forall_{B \in (s, A)^\omega}
        \exists_{B^* \in (\emptyset, \tilde{B})^\omega} \bigg( \max(s) < \min(B^*) \wedge \forall_{i \in B^*}
        \mbox{ \( A \) is \(C\)-accepting \( s\cup\{i\} \) }
        \bigg)
      \end{equation*}
    \end{Lem}
    \begin{Prf}
      Let \( B^*_0 = A\) and \(P(s) = \{ t_1, \ldots, t_n \} \). For \( 1 \le m \le n \) and \( B_{m} = B^*_{m-1} \) apply the previous lemma for \( (t_m, s \setminus t_m) \) obtaining \( B^*_m \). Then \( B^*_n \) has the required property.
    \end{Prf}

    \begin{Cor}
      If \( A \) is \(C\)-accepting \(\emptyset\) then there exists an increasing sequence \((a_n)_{n<\omega}\) such that
      \begin{equation*}
        \forall_{n < \omega} \mbox{ \( A \) is \(C\)-accepting \(\{ a_0, \ldots, a_n \}\)}.
      \end{equation*}
    \end{Cor}

    \begin{Cor}
      \label{cor2}
      If \( A \) is \(C\)-accepting \(\emptyset\) then there exists \( B \in [A]^\omega \) such that \( X_B \) has an unconditional basis.
    \end{Cor}
    \begin{Prf}
        By the previous corollary we have an increasing sequence \((a_n)_{n<\omega}\) such that
        \begin{equation*}
        \forall_{n < \omega} \mbox{ \( A \) is \(C\)-accepting \(\{ a_0, \ldots, a_n \}\)}.
        \end{equation*}
        Clearly \( (x_{a_n})_{n<\omega} \) is linearly independent. Let \( (b_n)_{n<\omega} \), \( (\varepsilon_n)_{n<\omega} \) be sequences of real numbers such that \( \sum\limits_{n = 0}^{\infty} b_n x_{a_n} \) is convergent and \( \varepsilon_n = \pm 1 \). For any \( n < \omega \) we have \( \| \varepsilon_0 b_0 x_{a_0} + \ldots \varepsilon_n b_n x_{a_n} \| \le C \| b_0 x_{a_0} + \ldots b_n x_{a_n} \| \) and consequently we also have \( \| \sum\limits_{n = 0}^{\infty} \varepsilon_n b_n x_{a_n} \| \le C \| \sum\limits_{n = 0}^{\infty} b_n x_{a_n} \| \), which shows that the sequence \( (x_{a_n})_{n < \omega} \) is unconditional. Now let \( (b_n)_{n < \omega} \) be such that \( \sum\limits_{n = 0}^{\infty} b_n x_{a_n} = 0 \). Then for any \( n < \omega \) we have \( \sum\limits_{k = 0}^{n-1} b_k x_{a_k} - b_n x_{a_n}  + \sum\limits_{k = n+1}^{\infty} b_k x_{a_k} = 0 \) and consequently \( b_n = -b_n \), ie \( b_n = 0 \). That completes the proof that \( (x_{a_n})_{n < \omega} \) is an unconditional basis for the space \( \Cl(\Span(\{ x_{a_n} \colon n \in \omega \})) \).
    \end{Prf}

  \section{The dichotomy}
  
    If any \( A \in (\emptyset, \omega)^\omega \) satisfying \( \tilde{A} = A \) is \(C\)-accepting \(\emptyset\) for some \( C > 1 \)  then it contains a \(C\)-dichotomous subset \( A^* \in (\emptyset,A)^\omega \) satisfying \( \tilde{A^*} = A^* \). Clearly \( A^* \) is also \( C \)-accepting \( \emptyset \) and by corollary \ref{cor2} we have a subspace with an unconditional basis. Otherwise let \( A_0 \in (\emptyset, \omega)^\omega \) be \(2\)-dichotomous and satisfying \( A_0 = \tilde{A_0} \). For each \(A_n\) we can find a \(n+3\)-dichotomous subset \( A_{n+1} \in (\emptyset, A_n)^\omega \) satisfying \( \tilde{A_{n+1}} = A_{n+1} \). Clearly \( A_{n+1} \) is \(n+3\)-rejecting \(\emptyset\) ie
    \begin{equation}
      \forall_{S, T \in (\emptyset, A_{n+1})^\omega} \unc(S, T) \ge n + 3
    \end{equation}
    
    Let us pick \( i_n \in A_n \) in such a way that the set \( B = \{ i_n \colon n \in \omega \} \) is linearly independent. Clearly \( B \) is almost contained in each \( A_n \), ie \( |B \ A_n| < \omega \).
    
    \begin{Lem}
      The set \(B\) defined as above satisfies
      \begin{equation}
        \forall_{n > 0} \forall_{S, T \in (\emptyset, \tilde{B})^\omega} \unc(S, T) \ge n + 2
      \end{equation}
    \end{Lem}
    \begin{Prf}
       Let \( S, T \in (\emptyset, \tilde{B})^\omega \) and let \( s_n = B \setminus A_n \). Then \( S,T \in (\emptyset, A_n \setminus \tilde{s_n})^\omega \). As \( \tilde{A_n} = A_n \) then also \( \tilde{A_n \setminus \tilde{s_n}} = A_n \setminus \tilde{s_n} \).
        
      It follows that
      \[
        \unc(S, T) \ge \unc(S \setminus \tilde{s_n}, T \setminus \tilde{s_n}) \ge n + 2.
      \]
    \end{Prf}

    \begin{Cor}
      The set \(B\) satisfies
      \begin{equation}
        \forall_{C \in (1;+\infty)} \forall_{S, T \in (\emptyset, \tilde{B})^\omega} \unc(S, T) \ge C
      \end{equation}
      and thus \( X_B \) is hereditarily indecomposable.
    \end{Cor}
  
  \begin{center}
    Competing interests: The authors declare none
  \end{center}


\begin{thebibliography}{}
   	\bibitem{Gow} W.T. Gowers, A new dichotomy for Banach spaces, Geom. Funct. Anal. 6 (1996) 1083-1093
   	\bibitem{Ell} E. Ellentuck, A new proof that analytic sets are Ramsey, J. Symbolic Logic 39 (1974) 163-165
    \bibitem{RF} T. Figiel, R. Frankiewicz, R. Komorowski, C. Ryll-Nardzewski, On hereditarily indecomposable Banach spaces, Annals of Pure and Applied Logic 126 (2004) 293-299
    \bibitem{GP} F. Galvin, K. Prikry, Borel sets and Rasey's theorem, J. Symbolic Logic 38 (1973), 193-198
    \bibitem{Maur} B. Maurey, A Note on Gowers' Dichotomy Theorem, Convex Geometric Analysis 34 (1998)
  \end{thebibliography}
\end{document}